\documentclass[12pt]{article}
\usepackage{diagrams}
\usepackage{graphicx}
\usepackage{amssymb}
\usepackage{amsmath}

\newtheorem{prop}{Proposition}
\newtheorem{thm}{Theorem}
\newtheorem{corol}{Corollary}
\begin{document}
\title{Abstract Projective Lines}
\date{}
\author{Anders Kock\\ \small University of Aarhus }
\maketitle
\section*{Introduction}

For $V$ a vector space over a field $k$, 
one has the Grassmannian manifold $P(V)$ consisting 
of 1-dimensional linear subspaces of $V$. If $V$ is 
$n+1$-dimensional, $P(V)$ is a copy on $n$-dimensional projective 
space. For $n\geq 2$, $P(V)$ has a rich combinatorial structure, in 
terms of incidence relations (essentially: the lattice of linear 
subspaces), in fact, this structure is so rich  that one can 
essentially  reconstruct $V$ from the combinatorial structure. 

But for $n=1$, this combinatorial structure (in the form of a lattice),  
is trivial; as expressed by
R.\ Baer, ``{\em A line 
\ldots has no geometrical structure,  if considered as an isolated or 
absolute phenomenon, since then it is nothing but a set of points 
with the number of points on the line as the only invariant}\ldots'', 
\cite{baer} p.\ 71.

However, it is our contention that a projective line has another kind 
of structure, making it possible to talk about a projective line as a 
set equipped with a certain structure, in such a way that 
isomorphisms (projectivities) between projective lines are bijective 
maps which preserve this structure.

The structure we describe (Section \ref{APL}) is that of a groupoid (i.e.\ a category 
where all arrows are invertible), and with  certain properties. 
The fact that the coordinate projective line $P(k^{2})$, 
more generally, a projective space of the form $P(V)$ (and also the 
projective plane in the classical synthetic sense) has such groupoid structure, was 
observed in \cite{CAPS}, and  further elaborated on in \cite{DL}; we 
shall recall the relevant notions and constructions from \cite{CAPS} 
in Section \ref{PV}, and a crucial observation from \cite{DL} in 
Section \ref{3trans}.
The present note may be seen as a completion of some of the efforts 
of these two papers. 

\section{Groupoid structure on $P(V)$}\label{PV}
Let $k$ be a field and let  $V$ a 2-dimensional vector space over $k$. We 
have a groupoid ${\bf L}(V)$, whose set of objects  is the set of $P(V)$ 
of 1-dimensional 
linear subspaces of $V$, and whose arrows are the linear isomorphisms 
between these. For $A\in P(V)$, the  linear isomorphisms $A\to A$ are 
in canonical bijective correspondence with the invertible scalars,
$${\bf L}(V)(A,A)=k^{*};$$ 
on the other hand, if $A$ and $B$ are {\em distinct} 1-dimensional 
subspaces, then the linear isomorphisms $A\to B$ are all of the form 
``projection from $A$ to $B$ in a certain unique {\em direction} $C$'',
 with $C\in P(V)$ and $C$ distinct from $A$ and $B$. (This 
also works in higher dimensions, cf. \cite{CAPS} and \cite{DL}; one just has to require that $C$ 
belongs to the 2-dimensional subspace spanned by $A$ and $B$.) This is in fact 
a bijective correspondence, so ${\bf L}(V)(A,B)$ is canonically 
identified with the set $P(V)\backslash \{A,B\}$. Here is a picture 
from \cite{CAPS}:

\begin{center}
{\includegraphics[
natheight=1.943200in,
natwidth=2.889300in,
height=1.9821in,
width=2.9334in
]%
{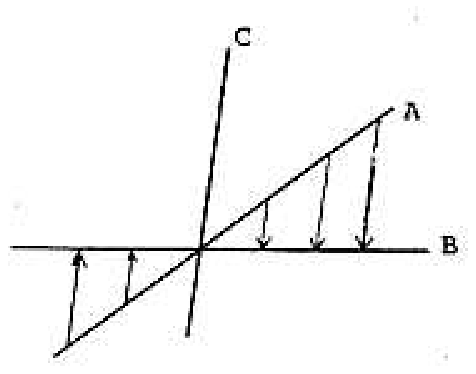}%
}
\end{center}
The linear isomorphism $A\to B$  thus described, we shall  denote 
$(C:A\to B)$. It is clear that the composite of $(C:A\to B)$ with $(C:B\to 
A)$ gives the identity map of $A$ (projecting forth and back in the 
same direction). Also it is clear that $(C:A\to B)$ composes with 
$(C:B\to D)$ to give $(C:A\to D)$. These equations will appear in the 
axiomatics for abstract projecive lines as the ``idempotency laws'', 
(\ref{idem1}) and (\ref{idem2}) below.

Also, it is clear that two linear isomorphisms from $A$ to $B$ differ 
by a scalar $\in k^{*}$; thus, for $A$ and $B$ distinct, and  
$(C:A\to B)$ and $(D:A\to B)$, there is 
a unique scalar $\mu \in k^{*}$ such that
\begin{equation}\mu . (C:A\to B) =(D:A\to B) = (C:A\to B).\mu 
.\label{star}\end{equation}
This scalar $\mu$ is (for $A,B,C,D$ mutually distinct) 
the classical cross-ratio $(A,B;C,D)$, cf.\ 
\cite{CAPS} (3) and \cite{DL} Theorem 1.5.3. (For $A,B,C$ distinct, 
and $D=C$, we have $(A,B;C,C) =1$.) Permuting the four entries 
(assumed distinct) 
will change the cross ratio according to well known formulae (see 
e.g.\ \cite{struik},\cite{KV}) which we 
shall make  explicit and take as axioms.

Thus, the groupoid ${\bf L}(V)$ which we in this way have associated to a 2-dimensional 
vector space $V$ over $k$ will be an example of an abstract projective 
line ${\bf L}$, in the sense of the next Section. 

\section{Abstract projective lines: axiomatics}\label{APL}

Let $k$ be a field. By a {\em $k$-groupoid}, 
we understand a groupoid ${\bf L}$ which is transitive (i.e.\ the hom 
set ${\bf L}(A,B)$ is non-empty, for any pair of objects $A,B$ in 
${\bf L}$), and such that all vertex groups ${\bf L}(A,A)$ 
are identified with the (commutative, multiplicative) group $k^{*}$  
of non-zero elements of the field $k$.  We assume
 that $k^{*}$ is {\em central} in ${\bf L}$ in the 
sense that for all $f:A\to B$ and $\lambda \in k^{*} = {\bf L}(A,A) = 
{\bf L}(B,B)$,
$\lambda .f = f\cdot \lambda$.
(We compose from the left to the right.) 

A $k$-functor between $k$-groupoids is a functor 
 which preserves 
$k^{*}$ in the evident sense.

We now define the notion of {\em abstract projective line} over $k$; it is 
to be a $k$-groupoid ${\bf L}$, equipped with the following kind of 
structure ($L$ denotes the set of objects of ${\bf L}$):

\noindent {\em for any two different objects $A,B\in L$, there is given a 
bijection between the  set ${\bf L}(A,B)$ and the set 
$L\backslash\{A,B\}$}
\noindent  satisfying some equational axioms: the {\em idempotence} laws 
(\ref{idem1}) and (\ref{idem2}), and the 
{\em permutation} laws (\ref{fourgroup}),\ldots ,(\ref{qqq}). To state these 
laws, we use, as in Section 1, the
notation:

 {\em if $C\in 
L\backslash\{A,B\}$, then the arrow $A\to B$ corresponding to it 
(under the assumed  bijection) by $(C:A\to B)$, or just by $C$, if $A$ and 
$B$ are clear from the context (say, from a diagram).}

Here are the first set of  equations that we assume (the 
``idempotence equations'') (we compose from left to right):
 Let $A,B,  F$ be mutually distinct, then

\begin{equation}(F:A \to B).(F:B\to A) = 1\in 
k^{*}\label{idem1}\end{equation}
and for $A,B,C,F$ mutually distinct
 \begin{equation}(F:A \to B).(F:B\to C) = (F:A\to 
C).\label{idem2}\end{equation}

The permutation laws which we state next are concerned with the 
crucial notion of {\em cross ratio}:
If $A,B,C,D$ are four distinct elements of $L$, we let $(A,B;C,D)$ 
be the unique scalar (element of $k^{*}$) such that
$$\begin{diagram}A&\rTo ^{C}& B\\
\dTo^{1}&&\dTo _{(A,B;C,D)}\\
A&\rTo _{D}&B
\end{diagram}$$
commutes; also, $(A,B;C,D)$ makes sense if $C=D$, and in this case 
equals $1\in k^{*}$. This scalar is called the {\em cross ratio} of 
the 4-tuple $A,B,C,D$.
\footnote{Convenience, as well as continuity, prompts us to define 
$(A,B;C,C)=1$ and $(A,B;C,B)=0$; this is consistent with  
determinant formulas for cross ratios in $P(k^{2})$ to be given 
later. In fact, one may consistently define $(A,B;C,D)$ whenever 
$A\neq D$ and $B\neq C$; $(A,A;C,D)=(A,B;C,C)=1$, and 
$(A,B;A,D)=(A,B;C,B)=0$.}

Since the elements of $L$ both appear as objects of ${\bf L}$ 
and as labels of arrows of ${\bf L}$, the four  entries (assumed 
distinct) in a 
cross ratio expression can be permuted freely by the 24 possible 
permutations of four letters. We assume the standard equation formulas for 
these permutation instances of a given cross ratio $\mu =(A,B;C,D)$; 
they give six values,
$$\mu , \frac{1}{\mu }, 1-\mu, \frac{1}{1-\mu} , \frac{\mu}{\mu -1}, \frac{\mu 
-1}{\mu},$$
see e.g. \cite{struik} p.\ 8 or \cite{KV} 0.2. The equations are
\begin{equation}(A,B;C,D)=(B,A;D,C)=(C,D;A,B)=(D,C;B,A),
    \label{fourgroup}\end{equation}
    and the following equations, where $\mu$ denotes $(A,B;C,D)$,
\begin{equation}(A,B;D,C)=\mu^{-1};	
	\label{inv}\end{equation}
\begin{equation}(A,C;B,D)= 1-\mu ; \quad (A,C;D,B)= \frac{1}{1-\mu}
    \label{transp}\end{equation}
    
	\begin{equation}(A,D;B,C)=\frac{\mu -1}{\mu};  \quad 
	    (A,D;C,B)=\frac{\mu}{\mu -1}.
	\label{qqq}\end{equation}
(This set of equations is not  independent.)
We had not needed to be so specific about these ``permutation equations'', 
since we shall only need the following 
consequence: if a map $\Phi :L \to L'$ preserves cross ratios of the form 
$(A,B;C,D)$ for some distinct $A,B,C,D$, then it also 
preserves any other cross ratio in which  the entries are 
$A,B,C,D$ in some other order.

\medskip

We have now stated what we mean by an abstract projective line ${\bf 
L}$. For (iso-){\em morphisms} (``projectivities'') between such: Let ${\bf L}$ and ${\bf L}'$ be abstract projective lines with 
object sets (underlying sets) $L$ and $L'$, respectively. By an
 {\em isomorphism} ${\bf L}\to {\bf L}'$ of projective lines, we 
understand a bijective map $\phi :L\to L'$ with the property that if 
we put 
$$\overline{\phi}(F:A\to B):=(\phi (F):\phi (A)\to \phi (B)),$$ 
(and $\overline{\phi}(\lambda )= \lambda$ for any scalar $\lambda \in 
k^{*}$), then
$\overline{\phi}$ commutes with composition, i.e.\ it defines a 
{\em functor} ${\bf L} \to {\bf L}'$ (preserving scalars, i.e.\ a 
$k$-functor).  The noticeable aspect of the 
category ${\mathcal L}$ of abstract projective lines, with (iso)morphisms as just 
defined, is that the ``underlying'' functor ${\bf L} \mapsto L$ (from ${\mathcal L}$ to the 
cate\-gory of sets) is a {\em faithful} functor, so that it makes sense 
to say whether a given function $L\to L'$ is a morphism (projectivity) or 
not. 

As always in such situations, it is convenient to use the same 
notation for the object itself, and its underlying set; so we 
henceforth do not have to 
distinguish notationally between ${\bf L}$ and $L$. 

 \medskip

Cross ratio was defined as a special case of composition; 
projectivities, in the sense defined here, commute with composition, since they 
are functors. Hence it is clear that projectivities preserve cross 
ratios.

\medskip
 In an (abstract) projective line ${\bf L}$, one may draw some diagrams that are 
meaningless in more general categories, like
the following square (whose commutativity actually can be {\em 
proved} on basis of the axiomatics):
\begin{equation}\begin{diagram}A&\rTo ^{C}&B\\
\dTo^{B}&&\dTo_{A}\\
C&\rTo _{-1}&C
\end{diagram}\label{minusone}\end{equation}
(where $A,B,C$ are three distinct points in ${\bf L}$). The commutativity 
of this diagram, for ${\bf L}=P(V)$, expresses an evident geometric 
fact that one sees by 
contemplating the figure (from \cite{CAPS}, p.\ 3):

\begin{center}{\includegraphics[
natheight=1.943200in,
natwidth=2.889300in,
height=1.9821in,
width=2.9334in
]%
{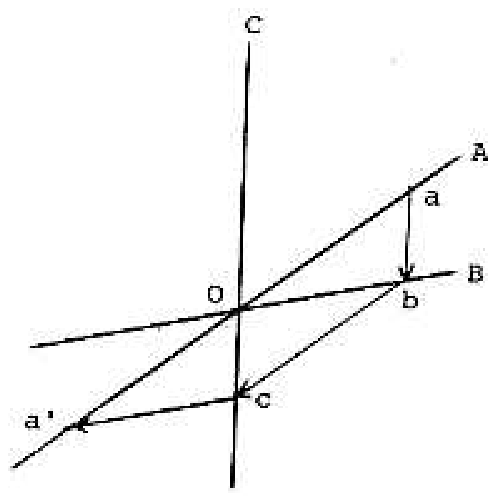}%
}
\end{center}
\noindent The existence of this diagram (\ref{minusone}) shows that ``cross ratios do 
not immediately encode all the geometry'' of projective lines; for, no cross ratio (except 1) can 
be concocted out of just three distinct points; four are needed.

\section{Three-transitivity}\label{3trans}

The ``Fundamental Theorem'' for projective lines derived from 
2-dimensional vector spaces is: for any two lists of three distinct 
points, there is a unique projectivity taking the points of the first 
list to the points of the second. This theorem, we shall prove holds 
for abstract projective lines.

Let ${\bf L}$ and ${\bf L'}$ be abstract projective lines over the 
field $k$.

\begin{thm}[Fundamental Thorem] Given three distinct points $A,B,C$  
in ${\bf L}$,  and given similarly $A', B', C'$ 
three distinct distinct points in ${\bf L}'$. 
Then there is a unique projectivity $\phi :{\bf L}\to {\bf L}'$ taking 
$A$ to $A'$, $B$ to $B'$ and $C$ 
to $C'$. 
\end{thm}
{\bf Proof.} For $D$ distinct from $A,B,C$, we put $\phi (D):=D'$, 
where $D' $ is the unique element in ${\bf L}'$ with $(A',B'; C',D') 
=(A,B;C,D)$; equivalently $D'$ is determined by the equation
$(C':A'\to B').(A,B;C,D) = (D':A'\to B')$.
By construction and the permutation equations, $\phi$ preserves cross 
ratios of any distinct 4-tuple, three of whose entries are $A,B,C$. 
Next, by the idempotence equations (\ref{idem1}) and (\ref{idem2}), 
$$(A,B;D;E)= (A,B;D,C).(A,B;C,E),$$
and similarly for the $A',\ldots, E'$. Each of the cross ratios on 
the right have three entries from the original set $A,B,C$, and so 
are preserved, hence so is the cross ratio on the left hand side,
$(A,B;D,E)$. So we conclude 
that any cross ratio,  two of whose entries are $A$ and $B$, is preserved.
Next,
$$(A,D;E,F)=(A,D;E,B).(A,D;B,F),$$
and similarly for the $A',\ldots ,F'$, so we conclude that any 
cross ratio with $A$ as one of its entries is preserved. 
Finally,
$$(D,E;F,G)=(D,E;F,A).(D,E;A,G),$$
and similarly for the $A', \ldots ,G'$, so we conclude that 
all cross ratios are preserved.

We have now described the bijection $\phi :{\bf L} \to {\bf L}'$, 
and proved that 
it preserves cross ratio of any four distinct points. To prove that 
it is a projectivity, in the sense defined, we need to argue that the 
corresponding $\overline{\phi}$ preserves composition of arrows. This 
is essentially an argument from \cite{DL} 2-4-4, which we make 
explicit:

\begin{prop}\label{DIL}If a bijection $\phi :{\bf L}\to {\bf L}'$ preserves 
cross ratio formation, 
then $\overline{\phi}$ preserves composition.
\end{prop}
It suffices to prove that commutative triangles go to commutative 
triangles. If the three vertices of the triangle agree,
these arrows are scalars $\in k^{*}$, and $\overline{\phi}$ preserves 
scalars. If  two, but not all three, vertices agree, one  arrow is a scalar, and 
commutativity of the triangle expresses that this scalar is  the 
cross ratio (or its inverse) of the four points that appear as the two 
vertices  and those two labels (likewise points in ${\bf L}$) that appear on 
the non-scalar arrows in the triangle; this is a an immediate 
consequence of the definition (\ref{star}), possibly combined with 
the idempotence law (\ref{idem1}). We conclude that composites of this 
form are likewise preserved by $\overline{\phi}$.   
Finally, we consider  the case where  the three vertices of 
the triangle are distinct, so the three arrows in the triangle  are 
of the form $(E:A\to B)$, $(F:B\to C)$, and $(G:A\to C)$ with $A, B, 
C$ distinct. 
Consider  $(E:A\to B).(F:B\to C). (G:C\to A)$, displayed as the top composite in 
the diagram 
$$\begin{diagram}
A&\rTo^{E}&B&\rTo ^{F}&C&\rTo ^{G}&A\\
\dTo^{1}&&\dTo_{1}&&\dTo_{1}&&\dTo_{(C,A;G,F)}\\
A&\rTo^{E}&B&\rTo^{F}&C&\rTo^{F}&A\\
\dTo^{1}&&\dTo^{1}&&&&\dTo_{(B,A;F,E)}\\
A&\rTo_{E}&B&&\rTo_{E}&&A.
\end{diagram}$$
All squares commute; the lower right hand rectangle commutes because 
of the idempotence law (\ref{idem2}) (the two $F$s combine into one). The lower 
composite is 1,  because of an idempotence law (ref{idem1}). So we conclude 
by (\ref{idem1}): $$(C,A;G,F).(B,A;F,E) =1 \mbox{  iff }(E:A\to B).(F:B\to C). (G:C\to 
A)=1.$$ Multiplying on the right by $G:A\to C$ (which is inverse to 
$G:C\to A$), we conclude
$$(C,A;G,F).(B,A;F,E) =1 \mbox{ iff } (E:A\to B).(F:B\to C)= (G:A\to 
C).$$
Thus commutativity of diagrams can be expressed in terms of cross 
ratio. Hence since cross ratio are preserved, the composite of 
$(E:A\to B)$ and $(F:B\to C)$  is preserved by 
$\overline{\phi}$. This proves the existence assertion of the 
Theorem. The uniqueness is clear, since a projectivity preserves 
cross ratios, so that we are forced to define $\phi (D)$ as the $D' 
\in {\bf L}'$ with $(A',B';C',D')=(A,B;C,D)$.

\section{${\bf L}=P(k^{2})$ as an abstract projective line}
The content of the present Section is mostly classical, but 
it emphasizes the category aspects of $P(k^{2})$. 
Non-zero vectors in $k^{2}$ are denoted $a=(a_{1},a_{2})$, $b=(b_{1},b_{2})$ 
etc.\, and the 1-dimensional linear subspaces of $k^{2}$ spanned by 
$a$ is denoted $A$; similarly, $b$ spans $B$, etc; $A$, $B$, \ldots are the 
points of the  set ${\bf L}= P(k^{2})$. We 
now have available the precious tool of {\em determinants} of 
$2\times 2$ matrices. We denote the determinant whose rows (or 
columns) are $a,b$ 
by the symbol $|a,b|$.

Given distinct $A$, $B$, and $C$, spanned by $a$, $b$, and $c$, respectively.
 We describe the linear map ``projection from 
$A$ to $B$ in the direction of $C$'' by describing its value  on $a\in 
A$; this value, being in $B$, is of the form $\lambda \cdot b$ for 
some unique scalar $\lambda \in k^{*}$, and an elementary calculation 
with linear equation systems (say, using Cramer's rule) gives that $\lambda =|c,a|/|c,b]$. Thus
\begin{equation}(C:A\to B)(a)= \frac{|c,a|}{|c,b|}\cdot 
b\label{basic}\end{equation}
is the basic formula. We can calculate the value of the composite
 $(C:A\to B).(D:B\to E)$
on $a\in A$; it takes $a\in A$ into
\begin{equation}\frac{|c,a|}{|c,b|}.\frac{|d,b|}{|d,e|}.e\in 
E.\label{compo}\end{equation}
In particular, if $E=A$, $a\in A$ goes into $(A,B;C,D). a$, where
$$(A,B;C,D):=\frac{|c,a|}{|c,b|}.\frac{|d,b|}{|d,a|}$$
$$=\frac{|a,c|.|b,d|}{|a,d|.|b,c |}$$ 
(using $|c,a|=-|a,c|$, and similarly for the other factors). This is 
the standard cross ratio $(A,B;C,D)$, and the standard permutation 
rules follow by known determinant calculations, as do the idempotency 
laws. So $P(k^{2})$ is indeed an abstract projective line, in 
our sense.

\medskip

In ${\bf L}=P(k^{2})$, 
we describe the points $A\in {\bf L}$ by homogeneous 
coordinates $[a_{1}:a_{2}]$, where $a$ is any vector spanning $A$. It is convenient to select three 
particular points
$V,H$, and $D$ (for ``vertical'', ``horizontal'', and ``diagonal'', 
respectively):
$$V=[0:1], H=[1:0], D=[1,1].$$
For any point $X$ distinct from $V$, there exists a unique $x\in k$ 
so that $X=[1:x]$. Thus, the $x\in k$ corresponding to $H$ and $D$ are $0$ 
and $1$, respectively. For $X$ distinct from $V$, the corresponding 
$x \in k$ may be calculated in terms of a cross ratio,
$$x=(V,H;D,X),$$
again by an easy calculation with determinants.
Thus ${\bf L}\backslash V$ has, by the chosen conventions,
 been put in 1-1 correspondence with the affine line $k$, so
$${\bf L}= \{ V\} + k;$$
$V$ is the ``point at infinity'' of the (``vertical'') copy 
$\{ (1,x) \mid  x\in k\}$ of the 
affine line $k$ inside $k^{2}$.

\medskip
The Fundamental Theorem then has the following 

\begin{corol}For every abstract projective line ${\bf L}$ over $k$, there 
exists  an isomorphism 
(=''projective equivalence'') 
with the projective line $P(k^{2})$. 
\end{corol}

(The isomorphism claimed is not unique, unless $k$ is the 2-element 
field.)
To prove the Corollary,
pick three distinct points $A,B,C$ in ${\bf L}$, and let $\phi$ be the 
unique projectivity (as asserted by the Theorem) to $P(k^{2})$ 
sending $A$ to $[1:0]$, $B$ to $[0:1]$, and $[C]$ to $[1:1]$.

\medskip

The isomorphism/projectivity $\phi$ described in this Corollary 
depends on the choice of $A,B,C$, and so is not canonical. However, 
it allows us to perform calculations in ${\bf L}$ using coordinates, 
in the form of such projective equivalence ${\bf L}\cong P(k^{2})$.

\medskip

Let us for instance prove  commutativity of (\ref{minusone}).
It suffices to prove that it  holds in ${\bf L}=P(k^{2})$. 
For, then it follows from the 
Fundamental Theorem that it also holds for three distinct points in 
an abstract projective line ${\bf L}$.
 
So consider points $A, B, C$ in $P(k^{2})$. Using (\ref{compo}), we see that the composite
$(C:A\to B).(A:B\to C)$ takes  $a\in A$ into
$$\frac{|c,a|}{|c,b|}.\frac{|a,b|}{|a,c|}.c,$$
and since $|c,a|=-|a,c|$,  two factors cancel except for the 
sign, and we are left with
$$-\frac{|a,b|}{|c,b|}.c=-\frac{|b,a|}{|b,c|}.c$$
which is the value of $-(B:A\to C)$ on $a$. (See \cite{DL} 1-4-2 for a more coordinate free 
proof.)

\medskip

To complete the comparison with the classical ``coordinate-'' 
projective line $P(k^{2})$,  we need to compare 
projectivities in our sense (functors) with classical projectivities, 
meaning maps $P(k^{2})\to P(k^{2})$ that are ``tracked'' by  linear 
automorphisms $k^{2}\to k^{2}$.

Let $f:k^{2}\to k^{2}$ be such linear automorphism. Then it defines a 
map $P(f):P(k^{2})\to P(k^{2})$ by $[a_{1}:a_{2}]\mapsto 
[f(a_{1}):f(a_{2})]$. We shall see  
 that this map preserves composition of arrows, hence is a functor; 
for, by (\ref{basic}),
$f(C:A\to B))$ takes $f(a)\in P(f)(A)$ to 
$$\frac{|f(c),f(a)|}{|f(c),f(b)|}.f(b) =
 \frac{|c,a|}{|c,b|}.f(b) \in P(f)(B)$$
(using the product rule for determinants and then cancelling the four 
occurrences of the determinant of $f$ that appear). The fact that 
composition is preserved is then a consequence of the formula 
(\ref{compo}).

On the other hand, every projectivity $\phi :P(k^{2}) \to P(k^{2})$ 
(in our sense) is of the form $P(f)$ for some linear automorphism 
$f:k^{2}\to k^{2}$ (which is in fact unique modulo $k^{*}$). Let 
$\phi (H)=A$, $\phi (V) =B$ and $\phi (D)=C$. Pick non-zero vectors $a\in A$, $b\in B$ 
and $c\in C$. The linear automorphism $f:k^{2}\to k^{2}$ with matrix
$$f=\left[ \begin{array}{cc}a_{1}&\lambda b_{1}\\
a_{2}&\lambda b_{2}\end{array}\right],$$
where $$\lambda :=-\frac{|c,a|}{|c,b|}$$
has the property that it takes $(1,0)$ to $a$, hence $P(f)$ takes 
$H$ to $A$; it takes $(0,1)$ to $\lambda b$, hence $P(f)$ takes $V$ 
to $B$; and finally, some calculation with Cramer's rule, say, shows 
that $f$ takes $(1,1)$ into a multiple of $c$, so $P(f)$ takes $D$ to 
$C$. Since $\phi$ and $P(f)$ both are projectivities, and they agree 
on $H,V$, and $D$, they agree everywhere, by the Fundamental Theorem. This proves that every 
projectivity $\phi :P(k^{2}) \to P(k^{2})$ (functor)
 is indeed tracked by a linear automorphism $k^{2}\to k^{2}$.

\medskip

{\bf Remark.} The projectivity  $\phi :P(k^{2}) \to P(k^{2})$ tracked by a linear 
automorphism $f:k^{2}\to k^{2}$ with matrix $[\alpha _{ij}]$ is  
also classically described as the {\em fractional linear transformation}
$$x\mapsto \frac{\alpha _{21}+ \alpha _{22}x}{\alpha _{11}+ \alpha 
_{12}x}.$$
This refers to the identification of $x\in k$ with $[1:x]\in P(k^{2})$.

\section{Structures on punctured projective lines}
\begin{prop}Given a projective line ${\bf L}$, and given $A\in {\bf L}$. 
Then ${\bf L}\backslash \{A\}$ carries a canonical structure of an affine line.
\end{prop}
{\bf Proof.} We first consider the case where ${\bf L}=P(k^{2})$, and where 
$A=V=[0:1])$. Now $P(k^{2}) \backslash V$ is identified with $k$ via 
$x\mapsto [1:x]$, and structure of affine line on $k$ gives by this 
identification a structure of affine line on $P(k^{2}) \backslash V$.
Using the Fundamental Theorem, the general result now follows
if we can to  prove that a projectivity $P(k^{2})\to P(k^{2})$ which 
fixes the point $V =[0:1])$ preserves affine combinations of the 
remaining points. 
A projectivity which fixes $V=[0:1]$ is tracked by a $2\times 2$ 
lower triangular matrix, or in terms of fractional linear 
transformations on $k$, by a function of the form
$x\mapsto (\alpha_{21}+\alpha _{22}x)/\alpha _{12}$, and this is an 
affine map $k\to k$, hence preserves affine structure (i.e.\ 
preserves linear combinations whose coefficient sum is 1).

\begin{prop}Given a projective line ${\bf L}$, and given $A, B\in 
{\bf L}$ 
with $A$, $B$ distinct. 
Then ${\bf L}\backslash \{A\}$ carries  a canonical structure of abstract 
vector line, with $B$ as $0$.
\end{prop}
{\bf Proof.} This is  analogous to the proof of the previous 
Proposition. The requirement that not only $V$, but also $H$ is 
preserved implies that the fractional linear transformation 
$(\alpha_{21}+\alpha _{22}x)/\alpha _{12}$ 
considered in the previous proof must have $\alpha_{21}=0$, and so is 
of the form $x\mapsto \alpha_{22}x/\alpha _{12}$, and hence is linear 
in $x$.

\medskip
Finally,
\begin{prop}Given a projective line ${\bf L}$, and given $A, B$, and 
$C\in {\bf L}$, mutually distinct. 
Then ${\bf L}\backslash \{A\}$ carries canonical structure of vector line 
with chosen basis vector, with $B$ as $0$ and $C$ as the chosen basis 
vector. 
\end{prop}
This last proposition is  a reformulation of the Fundamental 
Theorem. 

\section{The canonical bundles}
For each $A\in {\bf L}$, we have a canonical structure of affine line on 
${\bf L}\backslash \{A\}$. So over ${\bf L}$, we have a bundle ${\bf 
A}\to {\bf L}$
 of affine lines, whose fibre over $A\in {\bf L}$ is the 
affine line ${\bf L}\backslash \{A\}$. This bundle trivializes canonically over the 
covering ${\bf L}^{(3)}\to {\bf L}$, where ${\bf L}^{(3)}$ denotes the set of triples 
$A,B,C$ of mutually distinct points, and where the exhibited map is 
given by $(A,B,C)\mapsto A$.

The cocycle associated to this trivialization takes values in the 
group of affine automorphisms of $k$, which is a semidirect product 
of $(k,+)$ with $(k^{*},\cdot )$. If $h_{A,B,C}: {\bf L}\backslash 
\{A\} \to k$ is the restriction of the unique projectivity which 
takes $A$ to $[0:1]=V$, $B$ to $[1:0]=H$ and $C$ to $D=[1:1]$, then 
by construction of the affine structure on ${\bf L}\backslash 
\{A\}$, $h_{A,B,C}$ is an affine isomorphism. If $(A,B,C)$ and 
$(A,B',C')$ are two triples of distinct points (same $A$!), 
then the value of the desired cocycle on this pair is the affine isomorphism 
$h_{A,B,C}^{-1}. h_{A,B',C'}:k\to k$. Its value on $0$ is 
$h_{A,B',C'}(0)$, which is the cross ratio $(A',B';C',B)$, and 
similarly its value on 1 is the cross ratio $(A',B';C',C)$, so that 
the value of the cocycle on the pair $(A,B,C),(A,B',C')$ is
$$((A,B';C',B),(A,B';C',C)) \in k\ltimes k^{*}$$
(with $(A,B';C,B')$ by definition $=0$). 

\medskip

For each $A,B \in {\bf L}$, mutually distinct, we have a canonical vector 
line structure on ${\bf L}\backslash \{A\}$, with underlying affine line 
the one just described, and with $B$ as 0. 
So over ${\bf L}^{(2)}$ (= set of pairs of distinct points $A,B$ in ${\bf L}$), we have 
canonically a vector line bundle, whose fibre over $(A,B)$ is this 
vector bundle just described. It trivializes canonically over the 
covering
${\bf L}^{(3)}\to {\bf L}^{(2)}$ given by $(A,B,C)\mapsto (A,B)$ 
(here ${\bf L}^{(3)}$ is the set of  triples of distinct points in 
${\bf L}$).

The $k^{*}$-valued cocycle describing the associated principal bundle 
associates to $(A,B,C),(A,B,C')$ the cross ratio $(A,B;C,C')$.
The fact that this is a cocycle is the idempotence laws for cross 
ratio.

\medskip
Finally, over ${\bf L}^{(3)}$, we have a bundle of vector-lines-with-chosen 
basis-vector. This bundle is already itself trivial, since a vector 
line with chosen basis vector is uniquely isomorphic to $k$ with 
$1\in k$ as the chosen element.

\section{Stacks of projective lines}
The  notion of projective line, and of morphism 
(= isomorphism = projectivity) between such, as described here, is a (1-sorted) first 
order theory. This immediately implies that the notion of a {\em bundle} of 
projective lines over a space $M$ makes sense, and in fact, such bundles pull back 
along maps, and descend along  surjections, so projective line 
bundles form canonically a stack over the base category of sets, or, 
with suitable modifications, over the base category of 
spaces, say. Continuity, or other forms of cohesion,  will usually follow by the the fact 
that the 
constructions employed are canonical, as in \cite{SGM}, Section A.5). The study of 
bundles of 
of projective lines 
in the category of schemes, from \cite{KV}, was the input challenge 
for the present work, and I hope to push further into loc.\ cit.\ 
using the abstract-projective-line concepts.)
\medskip

\noindent{\bf Example.} Let $k$ denote the field of three elements 
${\mathbb Z}_{3}$. Every 4-element ${\bf L}$ set carries a unique structure of 
projective line over this $k$. We invite the reader to construct this 
structure (a groupoid with 4 objects, and  each hom-set a 2-element 
set); the composition laws follow from the idempotence equations; the 
cross ratio of the the four distinct points (in any order) is $-1$.

(Another argument: the group $PGL(2;{\mathbb Z})$ has 24 elements, 
which is also the number of permutations of a 4-element set, hence 
every permutation is a projectivity.)

\medskip

It follows that for any space $M$, and for any 4-fold covering $E\to 
M$, the bundle $E\to M$ is uniquely a bundle of projective lines over 
$k$. Clearly, such $E\to M$ need not have a section $M\to E$, so does 
not come about from a bundle of affine lines over $M$, by completing 
the fibres by points at infinity (the fibrewise infinity points would 
provide a cross section).




\begin{thebibliography}{99}
\bibitem{baer}R.\ Baer, Linear Algebra and Projective Geometry, 
Academic Press 1952.
\bibitem{DL}Y.\ Diers and J.\ Leroy, Cat\'{e}gories des points d'un espace 
projectif, Cahiers de Topologie et G\'{e}ometrie Diff\'{e}rentielles 
Cat\'{e}goriques 35 (1994), 2-28.
\bibitem{CAPS} A.\ Kock, The Category Aspect of Projective Space, 
Aarhus Universitet Matematisk Institut, Preprint Series 1974/75 No.\ 7.
\bibitem{SGM}A.\ Kock, Synthetic Geometry of Manifolds, Cambridge Tracts in 
Mathematics Vol.\ 180, 2010.
\bibitem{KV}J.\ Kock and I.\ Vainsencher, An Invitation to Quantum 
Cohomology, Progress in Mathematics, Vol.\ 
249, Birkh\"{a}user 2007.
\bibitem{struik}D.\ Struik, Analytic and Projective Geometry, 
Addison-Wesley Publ.\ Co.\ 1953.
\end{thebibliography}
\end{document}